\documentclass[twoside,a4paper,12pt,reqno] {amsart}
\usepackage{amsmath,mathrsfs}
\usepackage{epsfig}\usepackage{psfrag}\usepackage{subfigure}
\usepackage{tabularx}\usepackage{rotating}
\newtheorem{theorem}{Theorem}

\newtheorem{remark}{Remark}

\begin{document}
\title[The Banach-Saks Property of the Banach Product Spaces]
{The Banach-Saks Property of the Banach Product Spaces} 

\author{Zhenglu Jiang}
\address{Department of Mathematics, Zhongshan University, 
Guangzhou 510275, China}
\email{mcsjzl@mail.sysu.edu.cn}

\author{Xiaoyong Fu}
\address{Department of Mathematics, Zhongshan University, 
Guangzhou 510275, China
}
\email{mcsfxy@mail.sysu.edu.cn}


\baselineskip 16pt

\date{\today.}



\begin{abstract}
In this paper we first take a detail survey of the study of the Banach-Saks property of Banach spaces 
and then show the Banach-Saks property of the product spaces  
generated by a finite number of Banach spaces having the Banach-Saks property. 
A more general inequality 
for integrals of a class of composite functions is also given by using this property.  
\end{abstract}

\maketitle
 
\section{Introduction}
\label{intro} 
We begin by introducing the Banach-Saks property of a Banach space. 
A Banach space  is said to be of the Banach-Saks property if its arbitrary bounded sequence 
contains a subsequence whose arithmetic means of all the finite previous terms orderly comprise a  
sequence convergent  in the norm of the banach space. 
A Banach space  is said to be of the weak Banach-Saks property if its arbitrary weak convergent sequence 
contains a subsequence whose arithmetic means of all the finite previous terms orderly comprise a  
sequence convergent  in the norm of the Banach space. 
Such a subsequence is usually called the Banach-Saks sequence. 
Although any weakly convergent sequence is bounded, 
not all the bounded sequences are weakly compact. Therefore 
the Banach-Saks property is not equivalent to the weak Banach-Saks one in some Banach spaces. 
However, in uniformly convex spaces, 
the Banach-Saks property is equivalent to the weak Banach-Saks one.
It can be proved that $c_0$, $l_1$ and $L_1$ are of the weak Banach-Saks property (see \cite{ws})  
and that a Banach space has the Banach-Saks property  
if and only if it has the weak Banach-Saks property and it is reflexive.   
It can be also known that there are spaces without the weak Banach-Saks property  
such as $C[0,1]$ and $l_\infty$ (see \cite{nrf} or \cite{js}).  

As was shown by Castillo \cite{jmfcem}, 
there is a long history of studying the Banach-Saks property of a Banach space.  
Riesz \cite{fr} found in 1910 that any norm-bounded set in the space $L^p[0,1]$ of functions 
whose $p$-powers are integrable in the interval $[0,1]$ is 
weakly compact for any given $p\in (1,+\infty).$  
An analogical proposition relative to the strong convergence is obviously false. 
Instead, Banach and Saks \cite{bs} proved in 1930 that any norm-bounded 
sequence of functions of the space $L^p[0,1]$ for any given $p\in (1,+\infty)$ 
contains a subsequence whose arithmetic means of all the finite previous terms orderly comprise a  
sequence convergent  in the norm of the Banach space; meanwhile, 
a similar result in the space $C[0,1]$ of continuous functions on the interval $[0,1]$ 
was shown by Gillespie and Hurwitz \cite{gh} 
and independently by Schreier \cite{js}.  
Schreier \cite{js} also showed that a  weak convergent sequence of 
indicator functions on the interval $[0,1]$  
contains no Banach-Saks subsequences; thus, not all spaces have 
the weak Banach-Saks property. 
In the same year, using any convex combination instead of the arithmetic means, 
Zalcwasser \cite{zz} formed a convergent sequence from any given sequence in the space $C[0,1].$ 
After two years, Banach \cite{b32} proved that there is a closed linear manifold of the space $C[0,1]$ 
which is in one-to-one linear isometric correspondence with any given separable Banach space (see \cite{jac}).
It is then a natural idea to replace the arithmetic means with any convex combination 
for any Banach space, that is, 
if a sequence  in a Banach space is weakly convergent to some point then
there are convex combinations form a sequence which 
converges in the norm of the Banach space to this point. 
This is a result obtained by Mazur \cite{sm} in 1933 and its proof is just an application of 
the Hahn-Banach theorem. 
Kakutani \cite{sk} extended 
the above result of Banach and Saks to any uniformly convex space in 1939. 
 Later, Nishiura and Waterman \cite{nw} proved in 1963 that any space with the Banach-Saks property must 
be reflexive and Szlenk \cite{ws} found in 1965 that any weakly convergent (instead of bounded) 
sequence of the space $L[0,1]$ has a Banach-Saks subsequence. 
Baernstein II \cite{ab} showed in 1972 that not 
all reflexive spaces have the Banach-Saks property. 
Farnum \cite{nrf} showed in 1974 that any weak convergent sequence in 
 the space of all continuous complex-valued functions on a compact metric space with sup norm  
contains no Banach-Saks subsequences. This result is the same as that of Schreier \cite{js}. 

There are also many very interesting results about the passage of the Banach-Saks property. 
Partington \cite{jrp} proved in 1977 that the Banach-Saks property passes from a special Banach space $X$ to $l_p(X)$  
and that this result does not hold if one replaces the Banach-Saks property with the weakly Banach-Saks property. 
Early results on the behaviour of the Banach-Saks property under interpolation 
are due to Beauzamy  and Heinrich. In both cases they deal with the classical real interpolation. 
Beauzamy \cite{bb} showed in 1978 that the uniform convexity carries the Banach-Saks property   
 from one of two compatible Banach spaces to their interpolation spaces but  the Banach-Saks property does not 
pass from one of two compatible Banach spaces to their interpolation spaces (or vice versa). 
Beauzamy \cite{bbs} also gave the interrelations between the Banach-Saks property 
and the so-called P$_{\text 1}$ one applied to Banach spaces as well as to operators on a Banach space. 
Heinrich \cite{sh} gaves in 1980 a unified proof of several known results 
including the passage of the Banach-Saks property from a Banach space to another Banach space 
whose identity belongs to the former space. 
After that, there are some counterexamples given by Bourgain and Schachermayer 
in regard to passage of the Banach-Saks property from a Banach space $X$ to $L_p(X).$
Guerre \cite{sg} gave in 1980  an exposition of the example of Bourgain that 
$L^2([0,1];E)$ may fail to have the Banach-Saks property when $E$ has it. 
Using the Lions-Peetre interpolation method, Schachermayer \cite{ws81} constructed in 1981 an example,  
as a tree-like modification of the Baernstein space,  of a space $E$ with the Banach-Saks property  
such that $L_{[0,1]}^2(E)$ fails to have this property.  
Many properties relevant to the Banach-Saks property can be also found in 
 the literature such as the book of Beauzamy and Laprest\'{e} \cite{bl}.
It is worth mentioning that one does not yet know whether the dual of any space 
with the Banach-Saks property has this property or not. 

Calder\'{o}n \cite{apc} originally developed in 1964 a complex interpolation method to obtain intermediate spaces. 
Cwikel and Reisner \cite{cr} extended in 1984 the above Beauzamy result 
on the passage of the uniform convexity to the case of the complex interpolation spaces.
Vignati \cite{mv} obtained in 1987 a generalization of the above result due to Cwikel and Reisner,  
using the complex interpolation method of Calder\'{o}n. 

Besides the classical methods mentioned above, 
the so-called Ramsey methods are applied into the study of the Banach-Saks property  (see \cite{jmfc},  \cite{jd}). 
The Banach-Saks property is relevant to the extraction of subsequences 
and the Ramsey methods are very effective in the study of properties 
depending on the choice of subsequences. 
Thus it is natural to use these methods for the Banach-Saks property of a Banach space. 
These devices come from an original idea of Ramsey \cite{fpr} on the formal logic.  
Many relevant papers can be found in the literature (e.g., \cite{ab}, \cite{emk}). 
 The basic Ramsey theorem in the hierarchy of choice principles states that any infinite set $A$ 
has the following partition property: for every partition $\{X,Y\}$ of the two-element subsets of $A,$  
there is an infinite subset $A_0$ of $A$ such that 
all of the two-element subsets of $A_0$ are contained in 
only one of $X$ or $Y$ (see \cite{ab}, \cite{emk}). Thus, given a sequence and  
 a property  that all its subsequences may have, 
either there is a subsequence such that all its further subsequences have  this property, or there is a subsequence without this property. 
 More specifically speaking, a bounded sequence in a Banach space 
either contains a subsequence such that all its further subsequences are Banach-Saks sequences,  
or has a subsequence without any further Banach-Saks subsequence. 

It is also known that  the Banach-Saks property can be carried by the uniform convexity 
from a finite number of Banach spaces to their product space. 
It is hence natural to ask whether the Banach-Saks property passes 
from a finite number of arbitrary Banach spaces to their product space. 
Besides the above survey of the study of the Banach-Saks property of Banach spaces, 
this paper is to show the Banach-Saks property of the Banach product spaces 
generated by a finite number of Banach spaces having the Banach-Saks property. 
As an application of this property, a more general inequality 
for integrals of a class of composite functions is also given. 

\section{Theorems and Their Proofs}
\label{tip} 
We first recall the definition of the product space of 
 a finite number of Banach spaces. 
Let $m$ be any given natural number greater than $1$ and $B_j$ a Banach space 
with its element $x^{(j)}$ for  $j=1,2,\cdots,m.$ 
The product space $B$ is defined as the set of 
all the $n$-tuples $x=(x^{(1)},x^{(2)},\cdots,x^{(n)})$ with their norm having the property that 
${\parallel}x{\parallel}\to 0$ is equivalent to ${\parallel}x^{(j)}{\parallel}\to 0$ for  $j=1,2,\cdots,m.$ 
Generally, $B$ is denoted by $B=B_1\times B_2\times\cdots B_m.$ 
 It can be also found that $B$ is a Banach product space. 
It can be easily known from \cite{jac} and \cite{sk} that 
if $B$ is a uniformly convex product of all the finite numbe of $B_j$ 
and all the $B_j$ are uniformly convex Banach spaces, then 
$B$ is uniformly convex and so $B$ has the Banach-Saks property. 
It is also clear that the product of a finite number of uniformly convex Banach spaces 
is not in general uniformly convex. However,   
we can now prove 
\begin{theorem}
Let the norm of an element $x=(x^{(1)},x^{(2)},\cdots,x^{(n)})$ 
of the above product space $B=B_1\times B_2\times\cdots B_m$ be defined by 
${\parallel}x{\parallel}=N({\parallel}x^{(1)}{\parallel},{\parallel}x^{(2)}{\parallel},\cdots,{\parallel}x^{(m)}{\parallel})$ 
where $N(a_1,a_2,\cdots,a_m)$ is a nonnegative continuous function of the nonnegative variables $a_i$ 
with the property that 
${\parallel}x{\parallel}\to 0$ is equivalent to ${\parallel}x^{(j)}{\parallel}\to 0$ for  $j=1,2,\cdots,m,$  
increasing in each variable separately 
and satisfying 
\begin{equation}
N(ca_1,ca_2,\cdots,ca_m)=cN(a_1,a_2,\cdots,a_m)
\label{con01}
\end{equation}
for all $ c\geq 0$ and 
\begin{equation}
N(a_1+b_1,a_2+b_2,\cdots,a_m+b_m)\leq N(a_1,a_2,\cdots,a_m)+N(b_1,b_2,\cdots,b_m)
\label{con02}
\end{equation} 
for all $b_j\geq 0$ ($j=1,2,\cdots,m$).  
Assume that all the spaces $B_j$  ($j=1,2,\cdots,m$) have the Banach-Saks property.
Then the product space $B$ has the Banach-Saks property as well. 
\label{bsp}
\end{theorem}
\begin{proof}
To prove this theorem, it suffices to consider the case when $m=2.$ 
Let the product space $B$ be denoted by $B=B_1\times B_2$ and $\{x_i=(x_i^{(1)},x_i^{(2)})\}_{i=1}^{+\infty}$ 
a bounded sequence of $B,$ where $x_i^{(j)}\in B_j$ for all $i$ and $j.$ 
Thus $\{x_i^{(j)}\}_{i=1}^{+\infty}$ is a bounded sequence in $B_j$ for $j=1,2.$ 
Since $B_1$ has the Banach-Saks property, by Ramsey's theorem, there is a subsequence 
denoted by $\{x_{i_r}^{(1)}\}_{r=1}^{+\infty}$ 
all of whose further subsequences are the Banach-Saks sequences. 
Since $B_2$ has the Banach-Saks property, by Ramsey's theorem, 
there is a subsequence of this sequence $\{x_{i_r}^{(2)}\}_{i=1}^{+\infty}$ 
all of whose further subsequences are the Banach-Saks sequences.
Notice that the subsequence of $\{x_{i_r}^{(1)}\}_{r=1}^{+\infty}$ corresponding to this subsequence 
of $\{x_{i_r}^{(2)}\}_{r=1}^{+\infty}$ is still the Banach-Saks sequence. 
We may denote without loss of generality  
by $\{x_{i_r}^{(2)}\}_{r=1}^{+\infty}$ this subsequence 
of $\{x_{i_r}^{(2)}\}_{i=1}^{+\infty}.$ 
Then $\{x_{i_r}=(x_{i_r}^{(1)},x_{i_r}^{(2)})\}_{r=1}^{+\infty}$ is the Banach-Saks sequence 
in the product space $B.$ This completes our proof. 
\end{proof}
\begin{remark}
Theorem \ref{bsp} still holds if the Banach-Saks property 
is replaced by the weak Banach-Saks property.  
Jiang and Fu \cite{jf} recently 
showed the weak Banach-Saks property of $(L_\mu^p({\bf R}^n))^m$ for any fixed natural number $m$ 
using the two following techniques with only minor adjustments: 
one is given by Banach and Saks \cite{bs} for any fixed $p\in(1,+\infty),$ 
and another by Szlenk \cite{ws} in the case when $p=1.$ 
\end{remark}
\begin{remark}
In the case for any fixed $p\in(1,+\infty),$ 
there is another brief proof of the Banach-Saks property of $(L_\mu^p({\bf R}^n))^m.$ This proof is as follows.
Since $L_\mu^p({\bf R}^n)$ is a uniformly convex Banach space for any fixed $p\in(1,+\infty),$
 $(L_\mu^p({\bf R}^n))^m$ is uniformly convex for any fixed real number $p\in(1,+\infty)$ 
and any fixed natural number $m$ (see \cite{jac}).      
Thus, by Kakutani's theorem \cite{sk}, $(L_\mu^p({\bf R}^n))^m$ has the Banach-Saks property 
for any fixed $p\in(1,+\infty)$ 
and natural number $m.$    
\end{remark}
\begin{remark}
As an application of this theorem, we can show inequalities for 
integrals of a class of composite functions (see \cite{jf}, \cite{jipam}, \cite{jmaa}). 
There are also some other results relevant to the Banach-Saks property of Banach spaces 
in the recent literature (e.g., \cite{dssf}, \cite{cn}, \cite{mio85}, \cite{mio87}, \cite{ss}).  
\end{remark}

Moreover, using Theorem \ref{bsp},  we can give a more general inequality as follows: 
\begin{theorem}  
Let  ${\bf R}^n$ be the usual vector 
space of real $n$-tuples $x=(x_1,x_2,\cdots,x_n)$ and 
$B=B_1\times B_2\times\cdots B_m$ 
the Banach product space defined as in Theorem \ref{bsp},
where all the spaces $B_j$ ($j=1,2,\cdots,m$) 
are the Banach spaces of functions on ${\bf R}^n$ and have the  Banach-Saks property,  
 $m$ and $n$ be two positive integers. 
Suppose that a sequence $\{u_i\}_{i=1}^{+\infty}$ 
 weakly converges in the Banach product space
 $B$ to $u$   
 as $i$ goes to infinity, and that any norm-convergent 
sequence in $B_j$ ($j=1,2,\cdots,m$) is   
 convergent for almost every element in ${\bf R}^n.$
Assume that all the values of $u$ and $u_i$ ($i=1,2,3,\cdots$)  
belong to a convex set $K$ in ${\bf R}^m$  and that 
$f(w)$ is a nonnegative continuous convex function from $K$ to  ${\bf R}.$ 
 Then 
\begin{equation}
\mathop{\underline{\lim}}\limits_{i\rightarrow+\infty}
 \int_{\Omega}f(u_i)d\mu\geq\int_{\Omega}f(u)d\mu \label{ineq}
\end{equation}
for any measurable set $\Omega\subseteq {\bf R}^n,$ where 
$\mu$  is a nonnegative Lebesgue measure of ${\bf R}^n.$
\label{thineq}
\end{theorem}
\begin{proof}
Put $\alpha_i=\int_{\Omega}f(u_i)d\mu$ ($i=1,2,\cdots $) and 
$\alpha=\mathop{\underline{\lim}}\limits_{i\rightarrow +\infty}\int_{\Omega}f(u_i)d\mu$ 
for all $\Omega\subseteq {\bf R}^n.$  
Then there exists a subsequence of $\{\alpha_i\}_{i=1}^{+\infty}$ such that 
this subsequence, denoted without loss of generality 
by $\{\alpha_i\}_{i=1}^{+\infty},$ converges to $\alpha$ as $i\rightarrow +\infty.$  

Since $u_i$ converges weakly in $B$ to $u,$ 
by Theorem~\ref{bsp},  
it is easy to see that there exists a subsequence $\{u_{i_j}:j=1,2,\cdots\}$ 
such that $\frac{1}{k}\sum_{j=1}^ku_{i_j} \rightarrow u$ 
in $B$ as $k\rightarrow +\infty.$  
Thus there exists a subsequence of $\{\frac{1}{k}\sum_{j=1}^ku_{i_j}: k=1,2,\cdots\}$ 
such that this subsequence (also denoted without loss of generality by 
$\{\frac{1}{k}\sum_{j=1}^ku_{i_j}: k=1,2,\cdots\}$) satisfies that, as $k\rightarrow +\infty,$
\begin{equation}
\frac{1}{k}\sum_{j=1}^ku_{i_j} \rightarrow u  \hbox{  a.e. in  }{\bf R}^n.
\label{conv}
\end{equation} 
On the other hand, since all the values of $\{u_i\}_{i=1}^{+\infty}$ and $u$ 
belong to the convex set $K$  in ${\bf R}^m$ and 
$f(w)$ is a nonnegative continuous convex function from $K$ to  ${\bf R},$  we have
\begin{equation}
f(\frac{1}{k}\sum_{j=1}^ku_{i_j})\leq \frac{1}{k}\sum_{j=1}^kf(u_{i_j}). 
\label{conf1}
\end{equation}
By (\ref{conf1}) and Fatou's lemma, it follows that 
\begin{equation}
 \int_{\Omega}\mathop{\underline{\lim}}\limits_{k\rightarrow+\infty}
f(\frac{1}{k}\sum_{j=1}^ku_{i_j})d\mu
\leq\mathop{\underline{\lim}}\limits_{k\rightarrow+\infty}
\frac{1}{k}\sum_{j=1}^k\int_{\Omega}f(u_{i_j})d\mu. 
\label{conf2}
\end{equation}
Combining (\ref{conv}) and (\ref{conf2}), we can know that 
\begin{equation}
 \int_{\Omega}f(u)d\mu
\leq\mathop{\underline{\lim}}\limits_{k\rightarrow+\infty}
\frac{1}{k}\sum_{j=1}^k\int_{\Omega}f(u_{i_j})d\mu\equiv 
\mathop{\underline{\lim}}\limits_{k\rightarrow+\infty}
\frac{1}{k}\sum_{j=1}^k\alpha_{i_j}.
\label{conf3}
\end{equation}
Finally, by using the property of the convergence of 
$\alpha_i$ to $\alpha,$  
(\ref{conf3}) gives (\ref{ineq}). 
This completes our proof.
\end{proof}

\vskip0.2cm
\noindent{\small {\bf  Acknowledgement.} 
This work was supported by grants of NSFC 10271121 and 
joint grants of NSFC 10511120278/10611120371 and RFBR 04-02-39026. 
This work was also sponsored by SRF for ROCS, SEM.  
We would like to thank the referee of this paper 
for his/her valuable comments on this work.}
{\small 

}
\end {document}